\documentclass[12pt]{article}
\usepackage{amsmath,amsfonts,amssymb,rotating,colordvi}
\usepackage{color}
\usepackage[unicode=true,
 bookmarks=true,bookmarksnumbered=true,bookmarksopen=true,bookmarksopenlevel=2,
 breaklinks=false,pdfborder={0 0 1},backref=false,colorlinks=true]
 {hyperref}

\hypersetup{pdftitle={Melham's Conjecture on Odd Power Sums of Fibonacci Numbers},
 pdfauthor={Brian Y. Sun, Matthew H.Y. Xie, Arthur L. B. Yang},
 pdfsubject={ },
 pdfkeywords={Fibonacci numbers; Lucas numbers; Melham's conjecture; Fibonacci polynomials; Lucas polynomials },
 linkcolor=black,citecolor=black,urlcolor=blue,filecolor=blue,pdfpagelayout=OneColumn,pdfnewwindow=true,pdfstartview=XYZ,plainpages=false}

\def\qed{\nopagebreak\hfill{\rule{4pt}{7pt}}}
\def\proof{\noindent {\it{Proof.} \hskip 2pt}}

\newtheorem{theo}{Theorem}[section]

\newtheorem{lemm}[theo]{Lemma}

\newtheorem{coro}[theo]{Corollary}
\newtheorem{conj}[theo]{Conjecture}

\numberwithin{equation}{section}

\allowdisplaybreaks

\begin{document}

\begin{center}
{\large \bf  Melham's Conjecture on Odd Power Sums of Fibonacci Numbers}%
\end{center}

\begin{center}
 Brian Y. Sun$^{1}$, Matthew H.Y. Xie$^{2}$, Arthur L. B. Yang$^{3}$\\[6pt]

Center for Combinatorics, LPMC-TJKLC\\
Nankai University, Tianjin 300071, P. R. China

Email:  $^{1}${\tt brian@mail.nankai.edu.cn}, $^{2}${\tt xiehongye@mail.nankai.edu.cn}, $^{3}${\tt yang@nankai.edu.cn}
\end{center}

\noindent\textbf{Abstract.}
 Ozeki and Prodinger showed that the odd power sum of the first several consecutive Fibonacci numbers of even order is equal to a polynomial evaluated at certain Fibonacci number of odd order. We prove that this polynomial and its derivative both vanish at $1$, and will be an integer polynomial after multiplying it by a product of the first consecutive Lucas numbers of odd order. This presents an affirmative answer to a conjecture of Melham.

\noindent{\bf AMS Classification 2010:} 11B39, 05A19.


\section{Introduction}

Let $F_n$ denote the $n$-th Fibonacci number, and let $L_n$ denote the $n$-th Lucas number.
It is well known that the Fibonacci numbers and the Lucas numbers satisfy the
same recurrence relation
\begin{align*}
F_n=F_{n-1}+F_{n-2},\qquad L_n=L_{n-1}+L_{n-2}, \qquad \mbox{for $n\geq 2$},
\end{align*}
with different initial values $F_0=0$, $F_1=1$, $L_0=2$ and $L_1=1$. The main objective of this paper is to prove the
following conjecture, which was proposed by Melham \cite{melham2008}.

\begin{conj}[{\cite[Conjecture 2.1]{melham2008}}]\label{conj-2}
For any positive integers $n,m$, the sum
\begin{equation}\label{formula-Melham-F}
L_1L_3L_5\cdots L_{2m+1}\sum_{k=1}^n F_{2k}^{2m+1}
\end{equation}
can be expressed as
$(F_{2n+1}-1)^2P_{2m-1}(F_{2n+1})$, where $P_{2m-1}(x)$ is a polynomial of degree $2m-1$ with integer coefficients.
\end{conj}

Melham's conjecture was motivated by a result of Clary and Hemenway \cite{cla}, who obtained that
 \begin{align*}
\sum_{k=1}^n F_{2k}^3=\left\{
\begin{array}{ll}
F_n^2L_{n+1}^2F_{n-1}L_{n+2}/4, & \mbox{if $n$ is even;}\\[5pt]
L_n^2F_{n+1}^2L_{n-1}F_{n+2}/4, & \mbox{if $n$ is odd.}
\end{array}
\right.
\end{align*}
The key ingredient in their proof is the following identity
\begin{align}
4\sum_{k=1}^n F_{2k}^3=(F_{2n+1}-1)^2(F_{2n+1}+2).\label{clary-hemenway-F}
\end{align}
This inspired Melham to consider the
sum $\sum_{k=1}^n F_{2k}^{2m+1}$ for general positive integer $m$.
From experimental results for small $m$, Melham inferred that
this sum has the property as stated in Conjecture \ref{conj-2}.
An explicit expansion for $\sum_{k=1}^n F_{2k}^{2m+1}$, as a polynomial in $F_{2n+1}$,
was obtained by Ozeki \cite{Oze} and Prodinger \cite{prod2008} independently.
For any positive integer $m$, let
\begin{align}
S_{2m+1}(x)&=\sum_{s=0}^m x^{2s+1}\left(\sum_{k=s}^m \frac{1+2k}{1+2s}\binom{2m+1}{m-k}\binom{k+s}{2s}\frac{(-5)^{s-m}}{L_{2k+1}}\right)\nonumber\\
&\quad -\sum_{k=0}^m (-1)^{m-k}\binom{2m+1}{m-k}\frac{F_{2k+1}}{L_{2k+1}}5^{-m}\label{S(x)}.
\end{align}
The Ozeki-Prodinger formula can be stated as follows.

\begin{theo}[\cite{Oze,prod2008}] For any positive integers $n,m$, we have
\begin{align}\label{Oze-Formula}
\sum_{k=1}^n F_{2k}^{2m+1}&=S_{2m+1}(F_{2n+1}).
\end{align}
\end{theo}

We would like to mention that Prodinger \cite{prod2008} went much further and evaluated the following power sums:
\begin{align}
\sum_{k=0}^n F_{2k+\delta}^{2m+\varepsilon} \quad \mbox{and} \quad \sum_{k=0}^n L_{2k+\delta}^{2m+\varepsilon},
\qquad \mbox{for $\delta,\varepsilon\in \{0,1\}$}.
\end{align}
Prodinger's work was further generalized by Chu and Li \cite{chuli2011}, who obtained polynomial representation formulas for power sums of the extended Fibonacci-Lucas numbers.

It is natural to hope that the polynomial $S_{2m+1}(x)$, appearing in the Ozeki-Prodinger formula,
may serve as a candidate for solving Conjecture \ref{conj-2}, namely,
\begin{itemize}
\item[(i)] the polynomial $S_{2m+1}(x)$ has a factor $(x-1)^2$, and
\item[(ii)] the polynomial $L_1L_3L_5\cdots L_{2m+1}S_{2m+1}(x)$ has only integer coefficients.
\end{itemize}
Cooper and Wiemann \cite{Wiemann2004} already proved that the constant term of
the polynomial $L_1L_3L_5\cdots L_{2m+1}S_{2m+1}(x)$ is an integer. It should be mentioned that their work preceded that of Ozeki and Prodinger. In this paper, we shall prove (i) and (ii) and hence give an affirmative answer to Melham's conjecture.

Our proof is also motivated by a partial answer to Conjecture \ref{conj-2}, which was given by Wang and Zhang \cite{wz2012}. Their approach uses the Fibonacci polynomials and the Lucas polynomials. Recall that the $n$-th Fibonacci polynomial $F_{n}(x)$ and the $n$-th Lucas polynomial $L_{n}(x)$
can be given by
\begin{align}\label{fl-def}
F_{n}(x)=\frac{\alpha(x)^n-\beta(x)^n}{\alpha(x)-\beta(x)}\quad \mbox{ and }\quad L_{n}(x)=\alpha(x)^n+\beta(x)^n,
\end{align}
where
\begin{align*}
\alpha(x)={\frac{x+\sqrt{x^2+4}}{2}},\,\beta(x)={\frac{x-\sqrt{x^2+4}}{2}}.
\end{align*}
It is easy to verify that  $L_{-n}(x)=(-1)^n L_{n}(x).$ Moreover, these polynomials satisfy the following recurrence relation:
\begin{align*}
F_n(x)=xF_{n-1}(x)+F_{n-2}(x),\qquad L_n(x)=xL_{n-1}(x)+L_{n-2}(x),
\end{align*}
where $n$ can be any integer. 
Wang and Zhang \cite{wz2012} obtained the following result.

\begin{theo}[{\cite[Corollary 2]{wz2012}}]
For any positive integers $n,m$, the sum
\begin{equation*}
L_1(x)L_3(x) \cdots L_{2m+1}(x)
\sum_{j=1}^n F_{2j}^{2m+1}(x)
\end{equation*}
can be expressed as $(F_{2n+1}(x)-x)H_{2m}(x; F_{2n+1}(x))$, where $H_{2m}(x; y)$ is a polynomial in two variables $x$ and $y$ with integer coefficients and degree $2m$ of $y$.
\end{theo}

Letting $x=1$, the above result could be considered as a big progress toward proving Melham's conjecture, since
$F_n(1)$ and $L_n(1)$ are just the Fibonacci number $F_n$ and the Lucas number $L_n$.
However, as noted by Wang and Zhang \cite{wz2012}, their method could not give a complete answer to Conjecture \ref{conj-2}. Instead, they obtained the following result, which is another conjecture of Melham \cite[Conjecture 2.2]{melham2008}.

\begin{theo}[{\cite[Corollary 4]{wz2012}}]\label{wang-zhang}
For any positive integers $n,m$, the sum
\begin{equation*}
L_1L_3 \cdots L_{2m+1}
\sum_{j=1}^n L_{2j}^{2m+1}
\end{equation*}
can be expressed as $(L_{2n+1}-1)Q_{2m}(L_{2n+1})$, where $Q_{2m}(x)$ is an integer polynomial of degree $2m$.
\end{theo}

The remainder of this paper is organized as follows. We have mentioned above that our proof of Conjecture \ref{conj-2} is based on the Ozeki-Prodinger formula. In Section 2, we shall show that $S_{2m+1}(x)$ has only integer coefficients after multiplying it by $L_1L_3L_5\cdots L_{2m+1}$, where $S_{2m+1}(x)$ is given by \eqref{S(x)}. In Section 3, we shall prove that $S_{2m+1}(x)$ has a polynomial factor $(x-1)^2$ for any positive integer $m$, and thus completely prove Conjecture \ref{conj-2}. In Section 4, we shall present another proof of Theorem \ref{wang-zhang}, which is very similar in spirit to the proof of Conjecture \ref{conj-2}.

\section{The integer coefficients}

Throughout this section, assume that $m$ is a positive integer. Let $S_{2m+1}(x)$ be defined as in \eqref{S(x)}. The main result of this section is as follows.

\begin{theo}\label{thm-integer}
The polynomial $L_1L_3L_5\cdots L_{2m+1}S_{2m+1}(x)$ has only integer coefficients.
\end{theo}

Note that Cooper and Wiemann \cite{Wiemann2004} has already proved the constant term of $L_1L_3L_5\cdots L_{2m+1}S_{2m+1}(x)$ is an integer. To prove Theorem \ref{thm-integer}, in view of \eqref{S(x)}, it suffices to show that, for any $0\leq s\leq m$, the following product
\begin{align*}
\left(L_{2s+1}L_{2s+3}\cdots L_{2m+1}\right)\cdot \left(\sum_{k=s}^m \frac{1+2k}{1+2s}\binom{2m+1}{m-k}\binom{k+s}{2s}\frac{(-5)^{s-m}}{L_{2k+1}}\right)
\end{align*}
is an integer, or equivalently,
\begin{align*}
5^{m-s}\, |\, \left(L_{2s+1}L_{2s+3}\cdots L_{2m+1}\right)\cdot \left(\sum_{k=s}^m \frac{1+2k}{1+2s}\binom{2m+1}{m-k}\binom{k+s}{2s}\frac{1}{L_{2k+1}}\right).
\end{align*}
Motivated by the work of Wang and Zhang \cite{wz2012}, we turn to consider the divisibility
of its polynomial version, and obtain the following result.

\begin{theo} \label{thm-divisibility}  For any $0\leq s< m$, let
\begin{align}
R_1(x;m,s)&=L_{2s+1}(x)L_{2s+3}(x)\cdots L_{2m+1}(x)\label{R1}\\[5pt]
R_2(x;m,s)&=\sum_{k=s}^{m}{\frac{1+2k}{1+2s}}{\binom{2m+1}{m-k}}{\binom{k+s}{2s}}{\frac{1}{L_{2k+1}(x)}}.
\label{R2}
\end{align}
and  ${R(x;m,s)=R_1(x;m,s)R_2(x;m,s)}$. Then $(x^2+4)^{m-s}\, | \, R(x;m,s)$, precisely,
the quotient $R(x;m,s)/(x^2+4)^{m-s}$ is an integer polynomial.
\end{theo}

In the following we shall concentrate on the proof of Theorem \ref{thm-divisibility}.
Note that
\begin{align}\label{eqn-comb-coeff}
\frac{1+2k}{1+2s}\binom{k+s}{2s}=2\binom{k+s+1}{2s+1}-\binom{k+s}{2s}=2\binom{k+s}{2s+1}+\binom{k+s}{2s},
\end{align}
which implies that $R(x;m,s)$ is an integer polynomial. Therefore, to prove Theorem \ref{thm-divisibility}, it suffices to show that
\begin{align}\label{eqn-vanish}
R^{(j)}(x;m,s)|_{x=2I}=0, \mbox{ for $0\leq j\leq m-s-1$,}
\end{align}
where $I^2=-1$ and $R^{(j)}(x;m,s)$ denotes the $j$-th derivative of $R(x;m,s)$ with respect to $x$.
By the famous Leibniz formula, we have
\begin{align*}
R^{(j)}(x;m,s)=\sum_{i=0}^j \binom{j}{i} R_1^{(j-i)}(x;m,s)R_2^{(i)}(x;m,s).
\end{align*}
We find the following surprising result, from which \eqref{eqn-vanish} immediately follows.

\begin{lemm}\label{lem-main} For any $0\leq s< m$ and  $0\leq i\leq m-s-1$, we have
\begin{align}\label{eqn-vanish2}
R_2^{(i)}(x;m,s)|_{x=2I}=0.
\end{align}
\end{lemm}

To prove the above lemma, in view of \eqref{R2}, we need to compute the high-order
derivatives of the reciprocal of ${L_{2k+1}(x)}$. To this end, we will use a result due to Leslie \cite{diff}.

\begin{lemm}[\cite{diff}]\label{lemma4-1}
For any positive integer $i$, we have
\begin{align}\label{eq-recip}
\left(\frac{1}{f(x)}\right)^{(i)}=\sum_{j=1}^{i}(-1)^j{{i+1}\choose{j+1}}{\frac{1}{f(x)^{j+1}}}(f(x)^j)^{(i)},
\end{align}
where all the derivatives are assumed to exist.
\end{lemm}

To apply Lemma \ref{lemma4-1} to the high-order
derivatives of $\frac{1}{L_{2k+1}(x)}$, we further need to compute the powers of ${L_{2k+1}(x)}$. It is easy to show that
\begin{align}\label{eqn-lucaspower}
(L_{2k+1}(x))^{j}=\frac{1}{2}\sum_{l=0}^{j}\binom{j}{l}
(-1)^{l(2k+1)}L_{(j-2l)(2k+1)}(x).
\end{align}
In fact, for any $j\geq 0$, we have
\begin{align*}
(z-{z^{-1}})^j&=\sum_{l=0}^j(-1)^l\binom{j}{l}z^{j-2l},\\
(z-{z}^{-1})^j&=\sum_{l=0}^j(-1)^{j-l}\binom{j}{l}z^{2l-j}.
\end{align*}
Therefore,
\begin{align*}
(z-{z}^{-1})^j=\frac{1}{2}\sum_{l=0}^j\binom{j}{l}
\left((-1)^lz^{j-2l}+(-1)^{j-l}z^{2l-j}\right).
\end{align*}
From this \eqref{eqn-lucaspower} follows by
setting $z=\alpha(x)^{2k+1}$. 

Now to prove \eqref{eqn-vanish2}, in view of \eqref{eq-recip} and \eqref{eqn-lucaspower}, there remains to evaluate the high-order derivative of the Lucas polynomials at $x=2I$.

\begin{lemm}\label{Lx2I}
For any  integer $n$, we have
\begin{align*}
L_n(x)|_{x=2I}=2I^n
\end{align*}
and
\begin{align} \label{eq-der-val}
(L_n(x))^{(i)}|_{x=2I}=|n|(i-1)!I^{n-i}\binom{|n|+i-1}{2i-1},~~\text{for}~1\leq i\leq |n|.
\end{align}
\end{lemm}

\proof We first consider the case of $n>0$. It is easy to check the validity of the first identity by \eqref{fl-def}.
For the high-order derivatives, we shall use the following expansion of $L_n(x)$:
\begin{align*}
L_n(x)&=\sum_{k\geq 0} {\frac{n}{n-k}}{{n-k}\choose k}x^{n-2k}.
\end{align*}
Thus,
\begin{align*}
(L_n(x))^{(i)}&=\sum_{k\geq 0}{\frac{n }{n-k}}i!{{n-k}\choose k}\binom{n-2k}{i}x^{n-2k-i}.
\end{align*}
Let $h(n,i)=(L_n(x))^{(i)}|_{x=2I}$, namely,
\begin{align*}
h(n,i)&=\sum_{k\geq 0}{\frac{n }{n-k}}i!{{n-k}\choose k}\binom{n-2k}{i}(2I)^{n-2k-i}.
\end{align*}
It is clear that $h(n,n)=n!$. Applying the function \textbf{sumrecursion} of the Maple package \textbf{sumtools}, we get the recurrence
\begin{align*}
h(n,i)=\frac{2(2i+1)I}{(n-i)(i+n)}h(n,i+1).
\end{align*}
Repeatedly using this recurrence, we obtain that
\begin{align*}
h(n,i)&=\frac{n(n+i-1)!(i-1)!}{(n-i)!(2i-1)!}I^{n-i}=n(i-1)!I^{n-i}\binom{n+i-1}{2i-1}.
\end{align*}

For the case of $n<0$, the desired result immediately follows from 
the case of $n>0$ and the relation $L_{n}(x)=(-1)^n L_{-n}(x)$. 
Moreover, it is straightforward to check the case of $n=0$. This completes the proof.
\qed

To prove Lemma \ref{lem-main}, we also need the following result.

\begin{lemm} \label{thm-key-sum}
For any $0\leq s<m$, we have
\begin{align}\label{SUM}
\sum_{j=s}^m(-1)^j\binom{2m+1}{m-j}\binom{j+s}{2s}p(j)=0,
\end{align}
where $p(x)$ is a polynomial of degree less than $2m-2s+1.$
\end{lemm}

\proof Let $T(m,s)$ denote the left hand side of \eqref{SUM}. Since $\binom{j+s}{2s}=0$ for $j<s$,
we have
\begin{align*}
 T(m,s)&=\sum_{j=0}^m(-1)^j\binom{2m+1}{m-j}\binom{j+s}{2s}p(j)\\[5pt]
 &=\sum_{j=0}^m(-1)^{m-j}\binom{2m+1}{j}\binom{m-j+s}{2s}p(m-j).
\end{align*}
Noting that $\binom{m-j+s}{2s}=0$ for $j>m$, we further get
\begin{align*}
 T(m,s)&=(-1)^{m+1} \sum_{j=0}^{2m+1}(-1)^{(2m+1)-j}\binom{2m+1}{j} \binom{m-j+s}{2s}p(m-j).
\end{align*}
Now consider $\binom{m-j+s}{2s}p(m-j)$ as a polynomial in $j$, say
$$\binom{m-j+s}{2s}p(m-j)=\sum_{l\geq 0} a_l(m,s)j^l.$$
By the hypothesis condition, this is a polynomial of degree less than $2m+1$.
Thus, we have
\[
 T(m,s)=(-1)^{m+1} \sum_{l\geq 0}a_l(m,s)\left(\sum_{j=0}^{2m+1}(-1)^{(2m+1)-j}\binom{2m+1}{j}j^l\right)=0,
\]
which follows from the well known identity
\begin{align}\label{stirling-use}
\sum_{j=0}^{k}(-1)^{k-j}{k \choose j}j^n=0, \,\,for \,\, n< k.
\end{align}
This completes the proof.
\qed

It should be mentioned that \eqref{stirling-use} is closely related to the Stirling numbers of the second kind, see Stanley \cite[p. 34]{stanley1995}. This formula also plays an important role in the proof of Dixon's identity given by Guo \cite{guo2003}.

We now are able to prove Lemma \ref{lem-main}.

\noindent \textit{Proof of Lemma \ref{lem-main}.}
By \eqref{R2}, we have
\begin{align*}
R_2^{(i)}(x;m,s)|_{x=2I}&=\sum_{k=s}^{m}{\frac{1+2k}{1+2s}}{\binom{2m+1}{m-k}}{\binom{k+s}{2s}}\left(\frac{1}{L_{2k+1}(x)}\right)^{(i)}|_{x=2I}.
\end{align*}
If $i=0$, by Lemma \ref{Lx2I} and Lemma \ref{thm-key-sum}, then we have
\begin{align*}
R_2^{(i)}(x;m,s)|_{x=2I}&=\sum_{k=s}^{m}{\frac{1+2k}{1+2s}}{\binom{2m+1}{m-k}}{\binom{k+s}{2s}}\frac{1}{2 I^{2k+1}}\\
&=\frac{1}{2 I(1+2s)}\sum_{k=s}^{m}{(-1)^k}{\binom{2m+1}{m-k}}{\binom{k+s}{2s}} (1+2k)\\
&=0
\end{align*}
If $i>0$, by  Lemma \ref{lemma4-1} and  \eqref{eqn-lucaspower} , we get
\begin{align}
R_2^{(i)}(x;m,s)|_{x=2I}=&\sum_{k=s}^{m}{\frac{1+2k}{1+2s}}{\binom{2m+1}{m-k}}{\binom{k+s}{2s}}\nonumber\\[6pt]
&\quad \times \sum_{j=1}^{i}(-1)^j{{i+1}\choose{j+1}}{\frac{((L_{2k+1}(x))^j)^{(i)}|_{x=2I}}{(L_{2k+1}(x))^{j+1}|_{x=2I}}}\nonumber\\[7pt]
=&\sum_{k=s}^{m}{\frac{1+2k}{1+2s}}{\binom{2m+1}{m-k}}{\binom{k+s}{2s}}\nonumber\\[6pt]
&\quad \times \sum_{j=1}^{i}(-1)^j{{i+1}\choose{j+1}}{\frac{((L_{2k+1}(x))^j)^{(i)}|_{x=2I}}{2^{j+1}I^{(2k+1)(j+1)}}}. \nonumber
\end{align}
It follows from \eqref{eqn-lucaspower} that
\begin{align*}
{\frac{((L_{2k+1}(x))^j)^{(i)}|_{x=2I}}{2^{j+1}I^{(2k+1)(j+1)}}} &=\frac{\sum_{l=0}^{j}\binom{j}{l}
(-1)^{l(2k+1)}(L_{(j-2l)(2k+1)}(x))^{(i)}|_{x=2I}}{{2^{j+2}I^{(2k+1)(j+1)}}}.
\end{align*}
Using equation \eqref{eq-der-val}, we get
\begin{align}
{\frac{((L_{2k+1}(x))^j)^{(i)}|_{x=2I}}{2^{j+1}I^{(2k+1)(j+1)}}}
 =&\sum_{l=0}^{j}\frac{(-1)^k(2k+1)(i-1)!}{{{2^{j+2}I^{(i+1)}}}}|j-2l|\nonumber\\[5pt]
 &\quad \times \binom{j}{l}
\binom{|j-2l|(2k+1)+i-1}{2i-1}.\nonumber 
\end{align}
Combining the above identities yields that
\begin{align}
R_2^{(i)}(x;m,s)|_{x=2I}
=&\sum_{k=s}^{m}{\frac{1+2k}{1+2s}}{\binom{2m+1}{m-k}}{\binom{k+s}{2s}}\nonumber\\[6pt]
&\quad \times \sum_{j=1}^{i}(-1)^j{{i+1}\choose{j+1}}\frac{(-1)^k(2k+1)(i-1)!}{{{2^{j+2}I^{(i+1)}}}}\nonumber\\[6pt]
&\quad \times \sum_{l=0}^{j}|j-2l|\binom{j}{l}
\binom{|j-2l|(2k+1)+i-1}{2i-1}.\nonumber
\end{align}
Changing the order of summation, we obtain
\begin{align*}
R_2^{(i)}(x;m,s)|_{x=2I}
=&\sum_{j=1}^{i} \sum_{l=0}^{j}\frac{(-1)^j|j-2l|(i-1)!}{{{2^{j+2}I^{(i+1)}{(1+2s)}}}}{{i+1}\choose{j+1}}\binom{j}{l} \\[5pt]
&\times \sum_{k=s}^{m}(-1)^k{\binom{2m+1}{m-k}}{\binom{k+s}{2s}} p(k;j,l,i),
\end{align*}
where
$$p(k;j,l,i)=\binom{|j-2l|(2k+1)+i-1}{2i-1}(2k+1)^2.$$
Clearly, $p(k;j,l,i)$, as a polynomial in $k$, is of degree $2i+1<2m+1$.
From Lemma \ref{thm-key-sum} we deduce that
$$\sum_{k=s}^{m}(-1)^k{\binom{2m+1}{m-k}}{\binom{k+s}{2s}} p(k;j,l,i)= 0$$
for any $l,j$, and hence
$R_2^{(i)}(x;m,s)|_{x=2I}=0$, as desired. This completes the proof. \qed

We proceed to prove Theorem \ref{thm-divisibility} and Theorem \ref{thm-integer}.

\noindent \textit{Proofs of Theorem \ref{thm-divisibility} and Theorem \ref{thm-integer}.}
Lemma \ref{lem-main} implies \eqref{eqn-vanish}, and also implies that
\begin{align}\label{eqn-conj-vanish}
R^{(j)}(x;m,s)|_{x=-2I}=0, \mbox{ for $0\leq j\leq m-s-1$}.
\end{align}
Thus, $R(x;m,s)$ has a polynomial factor $(x^2+4)^{m-s}$.
It is well known that integer polynomials must factor into integer polynomial factors.
Therefore, we have $(x^2+4)^{m-s}\, \mid \, R(x;m,s)$, since $R(x;m,s)$ has only integer coefficients.
This completes the proof of Theorem \ref{thm-divisibility}. We further let $x=1$ and obtain that
\begin{align*}
5^{m-s}\, |\, \left(L_{1}L_{3}\cdots L_{2m+1}\right)\cdot \left(\sum_{k=s}^m \frac{1+2k}{1+2s}\binom{2m+1}{m-k}\binom{k+s}{2s}\frac{1}{L_{2k+1}}\right).
\end{align*}
This completes the proof of Theorem \ref{thm-integer}. \qed

\section{The polynomial factor}

In this section, we aim to prove that $S_{2m+1}(x)$ in \eqref{S(x)} has a polynomial factor $(x-1)^2$. Towards this end, it suffices to prove that
$S_{2m+1}(1)=0$ and $S_{2m+1}'(1)=0$, where $S_{2m+1}'(x)$ denotes the derivative of $S_{2m+1}(x)$ with respect to $x$.

Let us first note several identities concerning the Fibonacci polynomials and the Lucas polynomials, which will be used later.

\begin{lemm}\label{three-eq} For any nonnegative integer $k$, we have
\begin{align}
\sum_{i=0}^k \frac{x^{2i+1}}{2i+1}\binom{k+i}{2i}&=\frac{L_{2k+1}(x)}{2k+1};\label{Chu-Li-L}\\[5pt]
\sum_{i=0}^k \frac{(-1)^i(x^2+4)^i}{2i+1}\binom{k+i}{2i}&=(-1)^k\frac{F_{2k+1}(x)}{2k+1};\label{Chu-Li-F}\\[5pt]
\sum_{i=0}^k(-1)^i(x^2+4)^i\binom{k+i}{2i}&=(-1)^k\frac{L_{2k+1}(x)}{x}. \label{Chu-Li-O}
\end{align}
\end{lemm}

\proof We first prove \eqref{Chu-Li-L}. The Fibonacci polynomials $F_n(x)$ can be expanded as follows:
\begin{align}
F_n(x)=\sum_{i=0}^{\lfloor(n-1)/2\rfloor}\binom{n-i-1}{i}x^{n-2i-1}, \quad n\geq 0,\label{Expansion-Koshy}
\end{align}
see Koshy \cite[(37.2)]{kshoy}. By \eqref{Expansion-Koshy}, we immediately get that
\begin{align}
\sum_{i=0}^k\binom{k+i}{2i}x^{2i}=F_{2k+1}(x).\label{Expansion-F}
\end{align}
From \eqref{fl-def} it is easy to deduce that
\begin{align*}
L_n'(x)=nF_n(x), \quad n\geq 0.
\end{align*}
Combining the above two identities and replacing $x$ by $t$, we get that
\begin{align*}
\sum_{i=0}^k\binom{k+i}{2i} t^{2i}=\frac{L_{2k+1}'(t)}{2k+1}.
\end{align*}
Integrate both sides of this equation with respect to $t$ from $t=0$ to $t=x$.
For the left-hand side, we have
\begin{align*}
\int_0^x\sum_{i=0}^k\binom{k+i}{2i}t^{2i}dt=\sum_{i=0}^k\binom{k+i}{2i}\int_0^x t^{2i}dt=\sum_{i=0}^k \frac{\binom{k+i}{2i}x^{2i+1}}{2i+1},
\end{align*}
and for the right-hand side, we have
\begin{align*}
\int_0^x\frac{L_{2k+1}'(t)}{2k+1}dt=\frac{L_{2k+1}(x)}{2k+1}.
\end{align*}
This completes the proof of \eqref{Chu-Li-L}.

We proceed to prove \eqref{Chu-Li-F} and \eqref{Chu-Li-O}.
By \eqref{fl-def}, it is easy to verify that
\begin{align*}
F_{2k+1}(I\sqrt{x^2+4})=\frac{(-1)^k L_{2k+1}(x)}{x},
\end{align*}
where $I^2=-1$.
Substituting $I\sqrt{x^2+4}$ for $x$ in \eqref{Chu-Li-L} and \eqref{Expansion-F}, and then applying the above equality, we immediately get \eqref{Chu-Li-F} and \eqref{Chu-Li-O}.
\qed

If we set $x=1$ in Lemma \ref{three-eq}, then we get the following result.
It should be mentioned that the first two identities also appeared in \cite{chuli2011}.

\begin{coro}\label{three-nb}
For any nonnegative integer $k$, we have
\begin{align}
\sum_{i=0}^k \frac{\binom{k+i}{2i}}{2i+1}&=\frac{L_{2k+1}}{2k+1};\label{Chu-Li-LL}\\
\sum_{i=0}^k \frac{(-5)^i}{2i+1}\binom{k+i}{2i}&=(-1)^k\frac{F_{2k+1}}{2k+1};\label{Chu-Li-FF}\\
\sum_{i=0}^k {(-5)^i}\binom{k+i}{2i}&=(-1)^k{L_{2k+1}}. \label{lucasid}
\end{align}
\end{coro}

We are now able to give the main result of this section.

\begin{theo} \label{thm-factor} For any positive integer $m$, the polynomial $S_{2m+1}(x)$ has a polynomial factor $(x-1)^2$, namely, $S_{2m+1}(1)=0$ and $S_{2m+1}'(1)=0$.
\end{theo}

\proof Letting $x=1$ in \eqref{S(x)}, we get that
\begin{align*}
S_{2m+1}(1)&=\sum_{i=0}^m \sum_{k=i}^m \frac{1+2k}{1+2i}\binom{2m+1}{m-k}\binom{k+i}{2i}\frac{(-5)^{i-m}}{L_{2k+1}}\\
&\quad -\sum_{k=0}^m (-1)^{m-k}\binom{2m+1}{m-k}\frac{F_{2k+1}}{L_{2k+1}}5^{-m}\\
&=\sum_{k=0}^m \sum_{i=0}^k  \left(\frac{(-5)^i}{1+2i}\binom{k+i}{2i}\right)\binom{2m+1}{m-k}\frac{(1+2k)(-5)^{-m}}{L_{2k+1}}\nonumber\\
&\quad -\sum_{k=0}^m (-1)^{m-k}\binom{2m+1}{m-k}\frac{F_{2k+1}}{L_{2k+1}}5^{-m}\\
&=0,
\end{align*}
the last step by (\ref{Chu-Li-FF}). Moreover, it follows from \eqref{S(x)} that
\begin{align*}
S'_{2m+1}(1)&=\sum_{k=0}^m \sum_{i=0}^k  \left({(-5)^i}\binom{k+i}{2i}\right)\binom{2m+1}{m-k}\frac{(1+2k)(-5)^{-m}}{L_{2k+1}}\\
&=\sum_{k=0}^m \binom{2m+1}{m-k}{(1+2k)(-1)^k(-5)^{-m}}\\
&=(-5)^{-m}\sum_{k=0}^m \binom{2m+1}{m-k}{(1+2k)(-1)^k}\\
&=0,
\end{align*}
the second step by \eqref{lucasid} and
the last step by Lemma \ref{thm-key-sum}. This completes the proof. \qed

Finally, we give a proof of Conjecture \ref{conj-2}.

\noindent \textit{Proof of Conjecture \ref{conj-2}.}
By the Ozeki-Prodinger formula, we have
$$L_1L_3L_5\cdots L_{2m+1}\sum_{k=1}^n F_{2k}^{2m+1}=L_1L_3L_5\cdots L_{2m+1}S_{2m+1}(F_{2n+1}).$$
Let
$$P_{2m-1}(x)=L_1L_3L_5\cdots L_{2m+1}S_{2m+1}/(x-1)^2.$$
From Theorem \ref{thm-integer} and Theorem \ref{thm-factor} it follows that
$P_{2m-1}(x)$ is a polynomial of degree $2m-1$ with only integer coefficients. This completes the proof.
\qed

\section{A new proof of Theorem \ref{wang-zhang}}

In this section we shall present a new proof of Theorem \ref{wang-zhang} along the spirit of our proof of
Conjecture \ref{conj-2}. Our proof is based on the following formula due to Prodinger \cite{prod2008}, who showed that, for any nonnegative integers $n,m$,
\begin{align}\label{eqll-chu-li}
\sum_{j=1}^n L_{2j}^{2m+1}&=\sum_{s=0}^m L_{2n+1}^{2s+1}\sum_{k=s}^m \frac{1+2k}{1+2s}\binom{2m+1}{m-k}\binom{k+s}{2s}\frac{1}{L_{2k+1}}-4^m.
\end{align}

Similar to \eqref{S(x)}, let
\begin{align}\label{eqn-w}
W_{2m+1}(x)=\sum_{s=0}^m x^{2s+1}\sum_{k=s}^m \frac{1+2k}{1+2s}\binom{2m+1}{m-k}\binom{k+s}{2s}\frac{1}{L_{2k+1}}-4^m.
\end{align}

Now we can give a proof of Theorem \ref{wang-zhang}.

\noindent \textit{Proof of Theorem \ref{wang-zhang}.}
Let $Q_{2m}(x)=L_1L_3 \cdots L_{2m+1}W_{2m+1}(x)/(x-1)$.
In view of \eqref{eqn-comb-coeff}, we see that $L_1L_3 \cdots L_{2m+1}W_{2m+1}(x)$ is an integer polynomial.
It suffices to show that $W_{2m+1}(x)$ has a polynomial factor $(x-1)$, namely $W_{2m+1}(1)=0$.
Note that
\begin{align*}
W_{2m+1}(1)&=\sum_{s=0}^m \sum_{k=s}^m \frac{1+2k}{{1+2s}}\binom{2m+1}{m-k}{\binom{k+s}{2s}}\frac{1}{L_{2k+1}}-4^m\\
&=\sum_{k=0}^m \left({\sum_{s=0}^k \frac{1}{1+2s}}{\binom{k+s}{2s}}\right)\binom{2m+1}{m-k}\frac{1+2k}{L_{2k+1}}-4^m\\
&=\sum_{k=0}^m {\frac{L_{2k+1}}{1+2k}}\binom{2m+1}{m-k}\frac{1+2k}{L_{2k+1}}-4^m\\
&={\sum_{k=0}^m \binom{2m+1}{m-k}}-4^m=0,
\end{align*}
where the second equality follows from \eqref{Chu-Li-LL} and the last equality
follows from the binomial theorem, see also \cite{chuli2011}.
Thus, $Q_{2m}(x)$ must be an integer polynomial of degree $2m$, and moreover  by Prodinger's formula
\eqref{eqll-chu-li}, we have
\begin{align*}
L_1L_3 \cdots L_{2m+1}
\sum_{j=1}^n L_{2j}^{2m+1}=(L_{2n+1}-1)Q_{2m}(L_{2n+1}).
\end{align*}
This completes the proof. \qed

\noindent{\bf Acknowledgements.} This work was supported by the 973 Project, the PCSIRT Project of the Ministry of Education and the National Science Foundation of China.

\end{document}